\documentclass[10pt,a4paper]{article}

\usepackage{gastex}

\usepackage{amsmath,amssymb,amscd,latexsym}
\usepackage{stmaryrd}
\usepackage{graphics, psfig,graphicx}
\usepackage{eufrak}
\usepackage{multind}

\usepackage[all,dvips]{xy}

\usepackage{theorem}

\numberwithin{equation}{section}
\newcommand\refeq[1]{(\ref{#1})}

\newtheorem{def.notation}{D\'efinition--Notations}[section]
\newtheorem{defff}{Definition}[section]

\newtheorem{prop.def}{Proposition--D\'efinition}[section]

\newtheorem{propriete.def}{Propri\'et\'e--D\'efinition}[section]
\newtheorem{property}{Property}[section]

\newenvironment{proof}{\noindent{\textsf{\underline{Proof}: }}}{$\blacksquare$}
\newtheorem{lemma}{Lemma}[section]
\newtheorem{theorem}{Theorem}[section]
\newtheorem{theorem*}{Theorem}

\theoremstyle{break}
{\theorembodyfont{\rmfamily}}
{\theorembodyfont{\rmfamily}}
{\theorembodyfont{\rmfamily}\newtheorem{example}{Example}[section]}
{\theorembodyfont{\rmfamily}}
{\theorembodyfont{\rmfamily}}
{\theorembodyfont{\rmfamily}}
{\theorembodyfont{\rmfamily}\newtheorem{NTO}{Notation}[section]}
{\theorembodyfont{\rmfamily}\newtheorem{remark}{Remark}[section]}
{\theorembodyfont{\rmfamily}}

\theoremstyle{plain}
\newtheorem{remark.num}{\sous{Remark}: }[section]

\newcommand\ie{i.e. }
\newcommand\eg{e.g. }

\newcommand\N{\mathbb{N}}
\newcommand\R{\mathbb{R}}

\newcommand\C{\mathbb{C}}

\newcommand\dd{\text{d}}

\newcommand\lent{\llbracket }
\newcommand\rent{\rrbracket }
\newcommand\NO[1]{\ensuremath{\Arrowvert #1 \Arrowvert}}
\newcommand\lNO[1]{\ensuremath{\left\Arrowvert #1 \right\Arrowvert}}

\newcommand\sous[1]{\underline{#1}}

\newcommand\dsurd[2]{\frac{\partial #1}{\partial #2} }

\newcommand\mc[1]{\ensuremath{\mathcal{#1}} }

\newcommand\un{\mathbb{I}} 

\newcommand\lang{\left\langle }
\newcommand\rang{\right\rangle }

\newcommand\sol{\mc{K} }

\newcommand\foret{\mathbb{F}}
\newcommand\arbre{\mathbb{A}}
\newcommand\cop{\tilde{\Delta}}
\newcommand\tree{\mathbb{T}}
\newcommand\tens{\mc{T}}

\newcommand\init{\mc{I}}

\title{Butcher series and control theory}
\author{Dikanaina Harrivel\footnote{LAREMA, UMR 6093, Universit\'e d'Angers, France, \textbf{dika@tonton.univ-angers.fr}}}

\begin{document}

\renewcommand\cop{\varpi}

\maketitle

\begin{abstract}
We show how solutions of a non--linear differential equation can be written as sum indexed by planar trees: 
the \emph{Butcher series}. Then we use that property in order to control non--linear differential equation. We show that 
if the linearized system is controllable then the system itself is controllable if the nonlinear term is small enough and we 
express explicitly the control as a sum indexed by planar tree which each terms is obtained by minimization of 
a functional. 
\end{abstract}

\vspace{ 1cm} 

\renewcommand\thetheorem{}
\renewcommand\theequation{\arabic{equation}}
\renewcommand\theremark{\arabic{remark}}
\section*{Introduction}
Butcher series are sums indexed by planar trees introduced by J. C. Butcher \cite{Butcher} in order to study and classify
\cite{Butcher.complements} Runge Kutta methods in numerical analysis. Ch. Brouder noticed \cite{Brouder.arbre1} that the 
structure which underlies Butcher's calculations is the Hopf algebra of rooted trees, and this Hopf algebra is exactly the Hopf 
algebra defined by D. Kreimer in his paper about renormalization \cite{Kreimer}. 

Butcher series gives a precise description of the solutions of a non linear differential equation. In this paper we show that 
Butcher series provides a way to find explicitly a control for a non--linear differential system when the linearized system 
is itself controllable. We apply this method in order to study a very simple problem and its control. 
But Ch. Brouder noticed (see for instance \cite{Brouder.BIT}) that Butcher series can be used in a very large class of 
situation including non linear PDE's. There are other works on control theory based on similar perturbative expansion, 
see \eg the papers of Matthias Kawski \cite{Kawski} and Matthias Kawski--H\'ector J. Sussmann \cite{Kawski.Suss}, although the 
point of view differs from ours.\\

Let $n\in\N$, $A$ be a $n\times n$ matrix, $A\in\mc{M}_n(\R)$, $T>0$, $f\in L^2((0,T),\R^n)$ and $\lambda\in\R$. 
Then consider the following problem 
\begin{equation}\label{PF.intro}\tag{$\mc{P}_\lambda$}
\left\{
\begin{array}{l}
\displaystyle{x\in\sol:=\mc{C}^0([0,T],\R^n)\cap H^1([0,T],\R^n)}\\
\displaystyle{x'=Ax+f+\lambda F(x)}\\
\displaystyle{x(0)=x^0\in\R^n.}
\end{array}
\right.
\end{equation}
Here $F: \sol\longrightarrow L^2((0,T),\R^n)$ is such that $F(x)=\sum_{p\ge 2}F_p(x,\ldots,x)$ where for all $p\ge 2$, 
$F_p$ is a $p$--linear map $F_p:\sol^{\otimes p}\longrightarrow L^2((0,T),\R^n)$ such that the power series 
$\vert F\vert(z):=\sum_{p\ge 0}\NO{F_p}z^p$ converges for all $z\in\C$. Then problem \refeq{PF.intro} can be solved using 
Butcher series. Let us introduce them briefly. \\

Planar trees are rooted trees drawn into the plane with the root on the ground. The external vertices are called \emph{leaves} 
and the other \emph{internal vertices}, we denote by $\vert b\vert$ the number of internal vertices of a planar tree $b$. 
We say that a planar tree is non--degenerate if and only if each internal vertex has 
at least two childrens. Let denote by $\tree$ the set of non degenerate planar trees. Then the solution $x$ of \refeq{PF.intro} 
can be written as a sum over planar trees \cite{Butcher}, \cite{Butcher.complements} 
$$
x=\sum_{b\in\tree}\lambda^{\vert b\vert}x(b)
$$
where for all $b\in\tree$, the function $x(b)$ is obtained by solving the linear problem ($\mc{P}_0$) (\ie \refeq{PF.intro} for 
$\lambda=0$) for various $x^0$ and $f$. \\

In this paper, we investigate the following case: we suppose that $(x^0,f)$ is itself a kind of Butcher series \ie $u:=(x^0,f)$ 
writes 
$$
u=(x^0,f)=\sum_{b\in\tree}\lambda^{\vert b\vert}u(b)
$$
where $u(b)$ does not depend on $\lambda$ and where the sum converges in the $\R^n\times L^2((0,T),\R^m)$ topology. 
We define a map $\Phi\ast u:\tree\longrightarrow\sol$ such that the following holds 
\begin{theorem}\textnormal{\textbf{\!\!\!\ref{theorem.Butcher.controle}} }\sl{
If $\lambda$ is small enough then the solution $x$ of \refeq{PF.intro} with 
$(x^0,f)=u=\sum_{b\in\tree}\lambda^{\vert b\vert}u(b)$ writes 
$$
x=\sum_{b\in \tree}\lambda^{\vert b\vert}(\Phi\ast u)(b)
$$ 
where the sum converges in the $\sol$ topology. 
}
\end{theorem}
\begin{remark}
\begin{itemize}
\item We construct the map $\Phi\ast u$ using a coproduct $\cop$ on the set of planar trees; $\Phi\ast u$ can be seen as a 
convolution product (in the algebraic sense \cite{Milnor.Moore}) for $\varpi$. The coproduct $\cop$ defines a bialgebra structure 
on $\tree$ but \emph{not} a Hopf algebra structure. 
\item Butcher series can describe the solutions of a large class of non linear problem (see \eg \cite{Butcher.QFT}, 
\cite{ARBRE2}). Hence theorem \ref{theorem.Butcher.controle} can be generalized. 
\end{itemize}
\end{remark}

This theorem may be useful in order to control the problem \refeq{PF}. 
Let  $m\in\N$, suppose that the source $f$ of problem \refeq{PF} writes $f=Bv$ where $B\in\mc{M}_{m,n}(\R)$ is a $m\times n$ 
matrix \ie consider the following problem 
\begin{equation}\label{PF2.intro}\tag{$\mc{P}_\lambda$}
\left\{
\begin{array}{l}
\displaystyle{x\in\sol}\\
\displaystyle{x'=Ax+Bv+\lambda F(x)}\\
\displaystyle{x(0)=x^0\in\R^n}
\end{array}
\right.
\end{equation}
The question we are interested in is: given $x^0\in\R^n$, is there a function $v\in L^2((0,T),\R^m)$ such that the solution $x$ of 
\refeq{PF2.intro} satisfies $x(T)=0$ ?\\ 

If $\lambda=0$ the answer is well known. It suffices to consider the following adjoint 
problem 
\begin{equation}\label{duale.intro}\tag{$\mc{P}'$}
\left\{
\begin{array}{l}
\displaystyle{y\in\sol}\\
\displaystyle{y'=-A^*y}\\
\displaystyle{y(T)=y^T\in\R^n}
\end{array}
\right.
\end{equation}
and compute $\frac{\dd}{\dd t}\lang x,y\rang$ for $x$ and $y$ solution of \refeq{PF2.intro} with $\lambda=0$ and 
\refeq{duale.intro} respectively. Then integrating from $0$ to $T$, we get 
$$
\lang x(T),y^T\rang=\lang x^0,y(0)\rang+\int_0^T\dd t\lang Bv(t),y(t)\rang
$$
Hence $x(T)=0$ if and only if the right hand side vanishes for all $y^T\in\R^n$. This formulation has a variational 
interpretation and leads to optimal control theory. In this case we know that there exists $v\in L^2((0,T),\R^m)$ such 
that $x(T)=0$ if 
and only if $A$ and $B$ satisfies the Kalman condition \cite{Zuazua.Micu}, \cite{Lee.Markus}
\begin{equation}\label{Kalman.intro}\tag{K}
\text{rank}(B,AB,\ldots,A^{n-1}B)=n
\end{equation}

Now consider the case $\lambda\neq 0$. Following the same steps, we get that for all $x$ and $y$ solution of \refeq{PF2.intro} 
and \refeq{duale.intro} respectively, 
$$
\lang x(T),y^T\rang=\lang x^0,y(0)\rang+\int_0^T\dd t\lang Bv(t),y(t)\rang +\lambda\int_0^T\dd t\lang F(x(t)), y(t)\rang.
$$
The basic idea is to look for a control $v$ in the form of a sum indexed by planar tree 
$v=\sum_{b\in\tree}\lambda^{\vert b\vert}v(b)$. Then using theorem \ref{theorem.Butcher.controle}, we show that 
the last identity leads to 
\begin{multline}\label{tres.important}
\lang x(T),y^T\rang=\lang x^0,y(0)\rang+\int_0^T\dd t\lang Bv(\circ)(t),y(t)\rang\\
+\sum_{\substack{b\in\tree\\ \vert b\vert\neq 0}}\lambda^{\vert b\vert}\left(
\int_0^T\dd t\lang Bv(b)(t),y(t)\rang +\int_0^T\dd t\lang F(b)(t), y(t)\rang\right)
\end{multline}
where for all $b\in\tree$, $F(b)$ is defined using the $v(c)$ such that $\vert c\vert<\vert b\vert$. 
Notice that the right hand side of \refeq{tres.important} depends only on the family $(v(b))_{b\in\tree}$, there are no other 
unknown quantities. Hence we get a real condition on the controlability of \refeq{PF2.intro} by a function $v$ of the form 
$v=\sum_{b\in\tree}\lambda^{\vert b\vert}v(b)$. \\

In this paper we try to deal with the right hand side of \refeq{tres.important} by annihilating each term of the sum over planar 
trees. More precisely  we search for a family of function $(v(b))_{b\in\tree}$ of $L^2((0,T),\R^m)$ such that 
\begin{gather}
\lang x^0,y(0)\rang+\int_0^T\dd t\lang Bv(\circ)(t),y(t)\rang=0\label{un.intro}\\
\forall b\neq\circ\  ;\ 
\int_0^T\dd t\lang Bv(b)(t),y(t)\rang +\int_0^T\dd t\lang F(b)(t), y(t)\rang=0\label{deux.intro}
\end{gather}
for any solution $y$ of \refeq{duale.intro}. \\

Directly from \refeq{un.intro} we find out that a necessary condition is that the linear system is controllable, and in fact we 
show that it is more or less sufficient. 
\begin{remark}
Notice that this is not the only possibility to annihilate the right hand side of \refeq{tres.important}. Identity 
\refeq{tres.important} provides a very interesting way to deal with non--linear systems such that the linearized system is not 
controllable. 
\end{remark}

We show that if the Kalman condition \refeq{Kalman.intro} is satisfied, we can define a family $(v(b))_{b\in\tree}$ of 
elements of $L^2((0,T),\R^m)$ by minimizing a family of functionals $J(b):\R^n\longrightarrow\R$, $b\in\tree$ such that  
the following theorem holds. 
\begin{theorem}\textnormal{\textbf{\!\!\!\ref{controle.simple}} }\sl{
If $\lambda$ is small enough then the sum $v:=\sum_{b\in\tree}\lambda^{\vert b\vert}v(b)$ is well defined in $L^2((0,T),\R^m)$ and 
provides a control for problem \refeq{PF2} \ie the solution $x$ of \refeq{PF2} corresponding to $(x^0,v)$ satisfies $x(T)=0$.   
}
\end{theorem}
\begin{remark}
In fact we have a more precise result with an explicit condition for the convergence of the sum 
$\sum_{b\in\tree} \lambda^{\vert b\vert}v(b)$ (see theorem \ref{theorem.controle.dim.finie}). 
\end{remark}

In the first section we introduce planar trees and a coproduct on the set of non degenerate planar trees and then 
we prove theorem \ref{theorem.Butcher.controle}. In the second section, we describe the control $v(b)$ 
and we prove theorem \ref{controle.simple}

\renewcommand\thetheorem{\thesection.\arabic{theorem}}
\renewcommand\theequation{\thesection.\arabic{equation}}
\renewcommand\theremark{\thesection.\arabic{remark}}

\section{Planar Trees}\label{planar.trees}
\begin{defff}\label{def.arbres.plans}\sl{
A \emph{planar tree} is an oriented connected finite graphe without loop together with an embedding into the plane; we suppose  
that the graph has a particular node that no edge points to; this node is called the \emph{ root }Êof the tree. 
}
\end{defff}
\begin{remark}\label{remarque.rooted}
The set of planar trees differs from the set of \emph{rooted trees}. For instance the following planar trees are different 
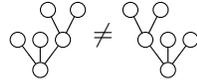
\begin{center}
\begin{picture}(0,3)(20,-2)
\gasset{Nw=2.0,Nh=2.0,Nmr=1.0,AHangle=0.0,AHLength=1.6}

\node(n1)(4,4){}
\node(n2)(8,4){}
\node(n3)(0,0){}
\node(n4)(3,0){}
\node(n5)(6,0){}
\node(n6)(3,-4){}
\drawedge(n1,n5){}
\drawedge(n2,n5){}
\drawedge(n3,n6){}
\drawedge(n4,n6){}
\drawedge(n5,n6){}

\put(10,-0.3){$\neq$}

\node(n11)(15,4){}
\node(n12)(19,4){}
\node(n13)(17,0){}
\node(n14)(20,0){}
\node(n15)(23,0){}
\node(n16)(20,-4){}
\drawedge(n11,n13){}
\drawedge(n12,n13){}
\drawedge(n13,n16){}
\drawedge(n14,n16){}
\drawedge(n15,n16){}

\end{picture}
\end{center}
although they represent the same rooted tree. 
\end{remark}
\begin{NTO}\label{notation.arbres.plans}\sl{
Let $b$ be a planar tree
\begin{enumerate}
\item The external vertices of $b$ are called \emph{leaves} and the other vertices are called \emph{internal vertices}. We will 
denote by $\NO{b}$ and $\vert b\vert$ respectively the number of leaves and internal vertices of $b$.
\item We denote by $\circ$ the planar tree without internal vertex and one leaf. 
\item A planar tree is non degenerate if each of its internal vertices has \emph{at least} $2$ childrens. In the following we 
consider only non--degenerate planar trees. We denote by $\tree$ the set of non--degenerate planar trees. 
\end{enumerate}
}
\end{NTO}
\begin{example}\label{exemple.nb}\sl{
The planar trees of remark \ref{remarque.rooted} are non degenerate and satisfy 
$\NO{b}=4$ and $\vert b\vert=2$. 
}
\end{example}

A very useful property is the recursive definition of planar trees: any planar tree which is not reduced to 
a single root can be obtained in a unique way by connecting the roots of $m$ trees to a new root. One can obtain all the planar 
trees by repeating this procedure. 

\begin{NTO}\label{algebre.tensorielle}\sl{
Let $V$ be a vector space then we denote by $\tens(V)$ the tensor algebra constructed over $V$ \ie 
$\tens(V):=\bigoplus_{p\ge 0}V^{\otimes p}$. Let us denote by $\arbre$ the $\R$--vector space spanned by $\tree$ and $\foret$ the 
tensor algebra over $\arbre$ \ie $\arbre:=\text{Vect}_\R\tree$ and $\foret:=\tens(\arbre)$. We denote by $\bullet$ the 
product on $\foret$ and by $\un:\R\longrightarrow\foret$ the unit of $\foret$.
}
\end{NTO}

\begin{defff}\label{def.B_+}\sl{
Consider $m\in\N$, $m\ge 2$ and $(b_1,\ldots,b_m)\in \tree^m$, then we denote by $B_+(b_1,\ldots,b_m)$ the planar tree obtained by 
connecting to a new root the roots of $b_1$ and $b_2$ and \dots and $b_m$. Hence $B_+$ can be considered as a linear map 
$B_+:\foret\longrightarrow\foret$.  
}
\end{defff}

\begin{property}\label{prop.recusrsion.arbre}\sl{
Let $b\in\tree$ be such that $b\neq\circ$ then there exists a unique $m\in\N$, $m\ge 2$ and $(b_1,\ldots,b_m)\in \tree^m$ such that 
$b=B_+(b_1,\ldots,b_m)$. We define $B_-:\arbre\longrightarrow\foret$ the linear map such that $B_-(\circ):=0$  and 
for $\forall b\in\tree$, $b\neq\circ$, $B_-(b):=b_1\bullet\cdots\bullet b_m$ where $b_1\bullet\cdots\bullet b_m$ is such that 
$b=B_+(b_1\bullet\cdots\bullet b_m)$. 
}
\end{property}

We now introduce a new operation on planar trees: the \emph{growing} operation and it's dual, the coproduct $\cop$. This operation 
leads to the perturbative expansion of solutions of nonlinear differential equations which we apply to control theory. 
\begin{defff}\label{excroissance}\sl{
Let $b\in\tree$ and $k$ denote the number of leaves of $b$. If $E=(E_1,\ldots,E_k)$ is a $k$-uplet $E\in \tree^k$ then 
we call the \emph{growing of $E$ on $b$} and denote by $E\propto b$ the planar tree obtained by replacing the $i$-th leaf of $b$ by
$E_i$ for all $i\in\lent 1,k\rent$. 
}
\end{defff}
\begin{example}\label{ex.excroissance}\sl{For instance we have 
\begin{center}
\begin{picture}(0,6)(50,-2)
\gasset{Nw=2.0,Nh=2.0,Nmr=1.0,AHangle=0.0,AHLength=1.6}
\node(n1)(5,4){}
\node(n2)(8,4){}
\node(n3)(11,4){}
\node(n4)(8,0){}
\node(n5)(12,0){}
\node(n6)(10,-4){}
\drawedge(n1,n4){}
\drawedge(n2,n4){}
\drawedge(n3,n4){}
\drawedge(n4,n6){}
\drawedge(n5,n6){}

\put(18,0){$=$}

\node(n7)(23,4){}
\node(n8)(26,4){}
\node(n9)(29,4){}
\node(n10)(26,0){}
\node(n11)(30,0){}
\node(n12)(28,-4){}
\drawedge(n7,n10){}
\drawedge(n8,n10){}
\drawedge(n9,n10){}
\drawedge(n10,n12){}
\drawedge(n11,n12){}
\put(32,0){$\propto$}
\node(n12)(37,0){}

\put(40,0){$=$}

\put(44,0){$($}
\node(n13)(47,2){}
\node(n14)(50,2){}
\node(n15)(53,2){}
\node(n16)(50,-2){}
\drawedge(n13,n16){}
\drawedge(n14,n16){}
\drawedge(n15,n16){}
\put(53,0){$,$}
\node(n17)(56,0){}
\put(58,0){$)$}

\put(60,0){$\propto$}
\node(n18)(65,2){}
\node(n19)(69,2){}
\node(n20)(67,-2){}
\drawedge(n18,n20){}
\drawedge(n19,n20){}

\put(72,0){$=$}

\put(76,0){$($}
\node(n21)(79,0){}
\put(81,0){$,$}
\node(n22)(84,0){}
\put(86,0){$,$}
\node(n23)(89,0){}
\put(91,0){$,$}
\node(n24)(94,0){}
\put(96,0){$)$}

\put(98,0){$\propto$}
\node(n25)(102,4){}
\node(n26)(105,4){}
\node(n27)(108,4){}
\node(n28)(105,0){}
\node(n29)(109,0){}
\node(n30)(107,-4){}
\drawedge(n25,n28){}
\drawedge(n26,n28){}
\drawedge(n27,n28){}
\drawedge(n28,n30){}
\drawedge(n29,n30){}

\end{picture}
\end{center}
}
\end{example}
\begin{defff}\label{coproduit}\sl{
Let $\cop:\foret\longrightarrow\foret\otimes\foret$ denote the morphism of algebras such that for all $b\in \tree$ 
$$
\cop(b):=\sum_{c\in \tree}\sum_{\substack{ E=(E_1,\ldots,E_{\NO{c}})\in \tree^{\NO{c}}\\ E\propto c=b}} 
(E_1\bullet\cdots\bullet E_{\NO{c}})\otimes c
$$
where $\foret\otimes\foret$ has the algebra structure inherited from $\foret$. 
}
\end{defff}
\begin{example}\label{exemple.cop}\sl{
For example, we have $\cop(1)=1\otimes 1$ because $\cop$ is an algebra morphism. By definition 
$\cop(\circ)=\circ\otimes\circ$ and using example \ref{ex.excroissance} we get 
\begin{center}
\begin{picture}(0,6)(50,-2)
\gasset{Nw=2.0,Nh=2.0,Nmr=1.0,AHangle=0.0,AHLength=1.6}
\put(-3,0){$\cop\biggl($}
\node(n1)(5,4){}
\node(n2)(8,4){}
\node(n3)(11,4){}
\node(n4)(8,0){}
\node(n5)(12,0){}
\node(n6)(10,-4){}
\drawedge(n1,n4){}
\drawedge(n2,n4){}
\drawedge(n3,n4){}
\drawedge(n4,n6){}
\drawedge(n5,n6){}
\put(14,0){$\biggl)$}

\put(18,0){$=$}

\node(n7)(23,4){}
\node(n8)(26,4){}
\node(n9)(29,4){}
\node(n10)(26,0){}
\node(n11)(30,0){}
\node(n12)(28,-4){}
\drawedge(n7,n10){}
\drawedge(n8,n10){}
\drawedge(n9,n10){}
\drawedge(n10,n12){}
\drawedge(n11,n12){}
\put(32,0){$\otimes$}
\node(n12)(37,0){}

\put(40,0){$+$}

\node(n13)(46,2){}
\node(n14)(49,2){}
\node(n15)(52,2){}
\node(n16)(49,-2){}
\drawedge(n13,n16){}
\drawedge(n14,n16){}
\drawedge(n15,n16){}
\put(54,0){$\bullet$}
\node(n17)(58,0){}
\put(60,0){$\otimes$}
\node(n18)(65,2){}
\node(n19)(69,2){}
\node(n20)(67,-2){}
\drawedge(n18,n20){}
\drawedge(n19,n20){}

\put(72,0){$+$}

\node(n21)(78,0){}
\put(80,0){$\bullet$}
\node(n22)(84,0){}
\put(86,0){$\bullet$}
\node(n23)(90,0){}
\put(92,0){$\bullet$}
\node(n24)(96,0){}
\put(98,0){$\otimes$}
\node(n25)(102,4){}
\node(n26)(105,4){}
\node(n27)(108,4){}
\node(n28)(105,0){}
\node(n29)(109,0){}
\node(n30)(107,-4){}
\drawedge(n25,n28){}
\drawedge(n26,n28){}
\drawedge(n27,n28){}
\drawedge(n28,n30){}
\drawedge(n29,n30){}

\end{picture}
\end{center}
}
\end{example}

\begin{remark}
We can prove that the coproduct $\cop$ leads to a bialgebra structure on $\foret$ which is \emph{not} a Hopf algebra structure. 
\end{remark}

\subsection{Application: Butcher series}\label{Butcher.series}
Let $T>0$ be a fixed positive real number, $n\in\N^*$ and $A\in\mc{M}_n(\R)$ be a $n\times n$ matrix. 
Then consider the following problem  
\begin{equation}\label{PF}\tag{$\mc{P}_\lambda$}
\left\{
\begin{array}{l}
\displaystyle{x\in\sol:=\mc{C}^0([0,T],\R^n)\cap H^1([0,T],\R^n)}\\
\displaystyle{x'=Ax+f+\lambda F(x)}\\
\displaystyle{x(0)=x^0\in\R^n}
\end{array}
\right.
\end{equation}
where $(x^0,f)$ belongs to $\init:=\R^n\times L^2((0,T),\R^n)$ and where $F$ is a map 
$F: \sol\longrightarrow L^2((0,T),\R^n)$ which satisfies the following hypothesis 
\makeatletter
\renewcommand{\theenumi}{\textbf{(H\arabic{enumi})}}
\makeatother
\begin{enumerate}
\item $F(x)=\sum_{p\ge 2}F_p(x,\ldots,x)$ where $F_p$ is a $p$--linear map 
$F_p:\sol^{\otimes p}\longrightarrow L^2((0,T),\R^n)$ for all $p\ge 2$. 
The power series $\vert F\vert$ defined by 
$
\vert F\vert(z):=\sum_{p\ge 0}\NO{F_p}z^p
$ 
converges for all $z\in\C$. 
\end{enumerate}
\makeatletter
\renewcommand{\theenumi}{\arabic{enumi}}
\makeatother

\begin{defff}\label{def.Phi.controle}\sl{
Let us define the family $(\Phi(b))_{b\in\tree}$ of $\NO{b}$--linear maps 
$\Phi(b):\init^{\otimes \NO{b}}\longrightarrow\sol$ recursively by setting 
$$
\Phi(\circ):(x^0,f)\in\init\longmapsto \text{solution of \refeq{PF} with }\lambda=0
$$ 
and for all $r\in\N^*$, $r\ge 2$ and for all $(b_1,\ldots,b_r)\in \tree^r$ 
\begin{equation}\label{def.rec}
\Phi\left(B_+(b_1,\ldots,b_r)\right):=\Phi(\circ)\bigl(0,\bullet\bigr) \circ F_r[\Phi(b_1)\otimes\cdots\otimes\Phi(b_r)]
\end{equation}
For all $b\in\tree$ we can see $\Phi(b)$ as a linear map $\Phi(b):\tens(\init)\longrightarrow\sol$. 
}
\end{defff}
Notice that since $\Phi(\circ)$ and $F_r$, $r\in\N^*$ are continuous, definition \refeq{def.rec} ensures that $\Phi(b)$ is 
a continuous $\NO{b}$--linear map $\Phi(b):\init^{\otimes \NO{b}}\longrightarrow\sol$.\\

We can show that for all $u=(x^0,f)\in\init$, if $\lambda$ is small enough then the power series 
$\sum_{b\in\tree}\Phi(b)(u)$ converges in $\sol$ and the sum is the solution of problem \refeq{PF}. 

Here we investigate a more general question: 
what happens if the initial condition $x^0\in\R^n$ or the source $f\in L^2((0,T),\R^n)$ depends on $\lambda$ or more specifically 
if they are infinite sum indexed by planar trees ? \\

Assume that $u:=(x^0,f)$ is a power series of the form 
$
u:=\sum_{b\in\tree} \lambda^{\vert b\vert} u(b)
$ 
where the family $(u(b))_{b\in\tree}$ of $\init$ is such that the power series $\vert u\vert$ defined by 
\begin{equation}\label{vert.u.vert}
\displaystyle{\vert u\vert:=\sum_{b\in\tree}\vert \lambda\vert^{\vert b\vert} \NO{u(b)}}
\end{equation}
converges. 

\begin{defff}\label{def.convolution.controle}\sl{
Let $\Phi\ast u$ denote\footnote{We can see $\Phi\ast u$ as the convolution product (in the algebraic sense \cite{Milnor.Moore}) 
of $\Phi$ and $u$ via the coproduct $\cop$.} the map $\Phi\ast u:\tree\longrightarrow\sol$ defined for all $b\in\tree$ by 
$$
(\Phi\ast u)(b):=\sum \Phi\left(b_{(2)}\right)\left(u(b_{(1)})\right)
$$
where we used the Sweedler notation $\cop(b)=\sum b_{(1)}\otimes b_{(2)}\in \foret\otimes \foret$. 
}
\end{defff}
\begin{theorem}\label{theorem.Butcher.controle}\sl{
If $u$ and $\lambda$ satisfy the following condition 
\begin{equation}\label{cond.Butcher.controle}
\vert \lambda\vert \vert u\vert^{-1}\ \vert F\vert\left(16 \NO{\Phi(\circ)}\vert u\vert\right)<1, 
\end{equation}
then the sum $\sum_{b\in \tree}\vert\lambda\vert^{\vert b\vert}\NO{(\Phi\ast u)(b)}$ converges and the sum 
$x=\sum_{b\in\tree} \lambda^{\vert b\vert}(\Phi\ast u)(b)$ is a solution of problem \refeq{PF}. 
}
\end{theorem}
\begin{remark}
As Ch. Brouder noticed \cite{Brouder.arbre1}, \cite{Brouder.BIT}, Butcher series can be used to solve a very large class of 
problem including PDEs (see \eg \cite{Butcher.QFT}) and we can prove a general version of theorem \ref{theorem.Butcher.controle}.  
\end{remark}

\begin{proof}(of theorem \ref{theorem.Butcher.controle})\\
Let us focus on the convergence of the sum. Looking at the definitions \ref{coproduit} and \ref{def.convolution.controle} 
we find out that it suffices to show that the sum 
\begin{equation}\label{but.dem1}
\sum_{\substack{b\in\tree\\ E\in\tree^{\NO{b}}}}\vert \lambda\vert ^{\vert b\vert+\vert E\vert}
\lNO{\Phi(b)\left(u(E)\right)}
\end{equation}
converges. For all $b\in\tree$, we denote by $N(b)$ the total number of vertices of $b$ \ie $N(b):=\NO{b}+\vert b\vert$.  
Let $N$ and $M$ belong to $\N^*$, then we have 
\begin{equation}\label{app.eq.Butcher.controle}
\sum_{\substack{b\in\tree\\ N(b)=N}}\vert \lambda\vert^{\vert b\vert}\sum_{\substack{E\in\tree^{\NO{b}}\\ N(E)\le M}}
\vert\lambda\vert^{\vert E\vert}\NO{\Phi(b)(u(E))}
\le 
\sum_{\substack{b\in\tree\\ N(b)=N}}\vert \lambda\vert^{\vert b\vert}\NO{\Phi(b)}
\sum_{\substack{E\in\tree^{\NO{b}}\\ N(E)\le M}}
\vert\lambda\vert^{\vert E\vert}\NO{u(E)}
\end{equation}
Since we assumed that the power series \refeq{vert.u.vert} converges, we have 
$$
\sum_{\substack{E\in\tree^{\NO{b}}\\ N(E)\le M}}
\vert\lambda\vert^{\vert E\vert}\NO{u(E)}\le \vert u\vert^{\NO{b}}
$$
So we can take the limit $M\to \infty$ in \refeq{app.eq.Butcher.controle} and get the following inequality 
\begin{equation}\label{truc}
\sum_{\substack{b\in\tree\\ N(b)=N}}\vert \lambda\vert^{\vert b\vert}\sum_{E\in\tree^{\NO{b}}}
\vert\lambda\vert^{\vert E\vert}\NO{\Phi(b)(u(E))}
\le 
\sum_{\substack{b\in\tree\\ N(b)=N}}\vert \lambda\vert^{\vert b\vert} \vert u\vert^{\NO{b}}\NO{\Phi(b)}
\end{equation}
Let us study $\NO{\Phi(b)}$. We denote by $I(b)$ the set of internal vertices of $b$ and 
for all internal vertex $i\in I(b)$ we denote by $r_b(i)\ge 2$ the number of childrens of $i$. Then we have the following lemma
\begin{lemma}\label{lemme.existence.controle}\sl{
For all $b\in \tree$ we have 
\begin{equation}\label{majoration.lemme.existence.controle}
\NO{\Phi(b)}\le \NO{\Phi(\circ)}^{ N(b)}\prod_{i\in I(b)}\NO{F_{r_b(i)}}
\end{equation}
}
\end{lemma}

\begin{proof}(of lemma \ref{lemme.existence.controle})\\
We will show \refeq{majoration.lemme.existence.controle} recursively on $N(b)$. If $N(b)=1$ then $b=\circ$, $I(b)=\emptyset$ and 
\refeq{majoration.lemme.existence.controle} is obvious. Fix $N\in\N^*$ and assume that \refeq{majoration.lemme.existence.controle} 
is true for all planar trees $b\in\tree$ such that $N(b)\le N$. Let $b\in\tree$ be such that 
$N(b)=N+1\ge 2$, then there exists $r\ge 2$ and $(b_1\ldots b_r)\in\tree$ such that $b=B_+(b_1\ldots b_r)$. Then definition 
\ref{def.Phi.controle} leads to 
$$
\lNO{\Phi(b)}\le \NO{\Phi(\circ)}\NO{F_r}\NO{\Phi(b_1)}\cdots \NO{\Phi(b_r)}
$$
and using \refeq{majoration.lemme.existence.controle} for the $b_i$'s we finally get 
\begin{equation}\label{fin.lemme}
\lNO{\Phi(b)}\le \NO{\Phi(\circ)}\NO{F_r}\prod_{j=1}^r \left(
\NO{\Phi(\circ)}^{\NO{b_j}}\prod_{i\in I(b_j)}\lNO{F_{r_{b_j}(i)}}
\right)
\end{equation}
But since $\vert b\vert=\vert b_1\vert+\cdots+\vert b_1\vert+1$, $\NO{b}=\NO{b_1}+\cdots+\NO{b_r}$ and since $I(b)$ is the disjoint 
union of the root of $b$ and the $I(b_j)$'s, we see that \refeq{fin.lemme} leads to 
\refeq{majoration.lemme.existence.controle}, which completes the proof. 
\end{proof}\\

Identity \refeq{truc} together with lemma \ref{lemme.existence.controle} leads to 
$$
\sum_{\substack{b\in\tree\\ N(b)=N}}\vert \lambda\vert^{\vert b\vert}\sum_{E\in\tree^{\NO{b}}}
\vert\lambda\vert^{\vert E\vert}\NO{\Phi(b)(u(E))}
\le 
\NO{\Phi(\circ)}^N\sum_{\substack{b\in\tree\\ N(b)=N}}\vert \lambda\vert^{\vert b\vert}
\vert u\vert^{\NO{b}}\prod_{i\in I(b)}\NO{F_{r_b(i)}}
$$
For all $b\in\tree$ such that $N(b)=N$, we have $N(b)=N=1+\sum_{i\in I(b)}r_b(i)$, so if we set $F_1=F_0=0$ we get 
\begin{multline}\label{dem.theo.existence.general.1}
\sum_{\substack{b\in\tree\\ N(b)=N}}\vert \lambda\vert^{\vert b\vert}\sum_{E\in\tree^{\NO{b}}}
\vert\lambda\vert^{\vert E\vert}\NO{\Phi(b)(u(E))}\\
\le 
\NO{\Phi(\circ)}^N \sum_{p=0}^N
\sum_{\substack{b\in \tree\\ \vert b\vert =p  ; \NO{b}=N-p}}
\sum_{\substack{(r_1,\ldots,r_p)\in\N^p\\  r_1+\cdots+r_p=N-1}}
\vert u\vert^{N-p}\vert \lambda\vert^p
\prod_{i=1}^p
 \NO{F_{r_i}}
\end{multline}
We know that the number of planar trees $c$ such that $N(c)=N$ is bounded by $16^N$ (see \cite{SEDGEWICK}), so the previous 
estimation leads to 
$$
\sum_{\substack{b\in\tree\\ N(b)=N}}\vert \lambda\vert^{\vert b\vert}\sum_{E\in\tree^{\NO{b}}}
\vert\lambda\vert^{\vert E\vert}\NO{\Phi(b)(u(E))}\le p_N
$$ 
where $p_N$ denotes the quantity 
\begin{equation}\label{dem.general.def.uN}
p_N:=
16^N\NO{\Phi(\circ)}^N\sum_{p=0}^N
\sum_{\substack{(r_1,\ldots,r_p)\in\N^p\\  r_1+\cdots+r_p=N-1}}
\vert u\vert^{N-p}
\vert \lambda\vert^p
\prod_{i=1}^p\NO{F_{r_i}}.
\end{equation}

Consider the formal power series $P:=\sum_{N\ge 1}p_N X^N\in \R[[X]]$.  
Using definition \refeq{dem.general.def.uN}, inverting the sum over $N$ and $p$ and considering that $F_0=F_1=0$ we get 
\begin{equation*}
P=
16X\NO{\Phi(\circ)}\vert u\vert\sum_{p\ge 1}\sum_{N\ge p}\sum_{\substack{(r_1,\ldots,r_p)\in\N^p\\  r_1+\cdots+r_p=N-1}} 
\left(
16X\NO{\Phi(\circ)}\vert u\vert
\right)^{N-1}\vert u\vert^{-p}\vert \lambda\vert^p
\prod_{i=1}^p\NO{F_{r_i}}
\end{equation*}
We recognize in the sum over $N$ the expression of 
$$
\left[\vert \lambda\vert \vert u\vert^{-1}\sum_{r\ge 0} \NO{F_r}\left(
16 X\NO{\Phi(\circ)}\vert u\vert
\right)^r\right]^p=\left[\vert \lambda\vert \vert u\vert^{-1}\ \vert F\vert\left(16 X\NO{\Phi(\circ)}\vert u\vert\right)\right]^p
$$ 
which tends to $0$ in the usual topology of $\R[[X]]$ since $F_0=F_1=0$. Hence we finally get the following identity 
\begin{equation*}
P=\frac{16X\NO{\Phi(\circ)}\vert u\vert}
{1-\vert \lambda\vert \vert u\vert^{-1}\ \vert F\vert\left(16 X\NO{\Phi(\circ)}\vert u\vert\right)}
\end{equation*}
So if condition \refeq{cond.Butcher.controle} is satisfied, then the radius of convergence of $P$ is larger than $1$ and the 
power series $\sum_N p_N$ converges which shows that the sum \refeq{but.dem1} converges. \\

Let us show that the sum $x=\sum_{b\in\tree}\lambda^{\vert b\vert}(\Phi\ast u)(b)$ is a solution of problem \refeq{PF}. 
We have shown that we have 
\begin{equation*}
\begin{aligned}
x=& \lim_{N\to\infty}
\left\{\lim_{M\to\infty}
\sum_{\substack{b\in\tree\\ \N(b)\le N}}\lambda^{\vert b\vert}\Phi(b)
\left(
\sum_{\substack{(E_1\ldots E_{\NO{b}})\in\tree^{\NO{b}}\\ N(E_j)\le M}}
\left[\lambda^{\vert E_1\vert}u(E_1)\right]\otimes\cdots\otimes \left[\lambda^{\vert E_{\NO{b}}\vert} u(E_{\NO{b}})\right]
\right)
\right\}\\
=&\lim_{N\to\infty}
\left\{\sum_{\substack{b\in\tree\\ \N(b)\le N}}\lambda^{\vert b\vert}\Phi(b) \left(u^{\otimes \NO{b}}\right)
\right\}
\end{aligned}
\end{equation*}
since $\sum_b \vert \lambda\vert^{\vert b\vert}\NO{u(b)}$ converges. Using definition\ref{def.Phi.controle} of $\Phi(b)$ we get 
that $\Phi(\circ)(u)=u_1=x^0$ and for all $b\neq \circ$, $\Phi(b)(u^{\otimes \NO{b}})(0)=0$, hence $x(0)=x^0$. Moreover for all 
$N\in\N^*$, we have 
\begin{multline*}
\dsurd{}{t}\left\{\sum_{\substack{b\in\tree\\ \N(b)\le N}}
\lambda^{\vert b\vert}\Phi(b) \left(u^{\otimes \NO{b}}\right)
\right\}= \sum_{\substack{b\in\tree\\ \N(b)\le N}}
\lambda^{\vert b\vert}A\Phi(b) \left(u^{\otimes \NO{b}}\right) \\ + 
\sum_{r\ge 2}\sum_{\substack{(b_1\ldots b_r)\in\tree^r\\ N(b_1)+\cdots+N(b_r)+1\le N}}
\lambda^{\vert b_1\vert+\cdots \vert b_r\vert +1}F_r\left(\Phi(b_1)(u^{\otimes \NO{b_1}})\otimes \cdots \otimes 
\Phi(b_r)(u^{\otimes \NO{b_r}})\right)
\end{multline*}
Then since 
$$
\sum_{\substack{(b_1\ldots b_r)\in\tree^r\\ N(b_1)+\cdots+N(b_r)+1\le N}}
\lambda^{\vert b_1\vert+\cdots \vert b_r\vert}F_r\left(\Phi(b_1)(u^{\otimes \NO{b_1}})\otimes \cdots \otimes 
\Phi(b_r)(u^{\otimes \NO{b_r}})\right)\substack{\longrightarrow\\ N\to\infty} \hspace{0,2cm}F_r(x^{\otimes r}),
$$ 
we get that $x$ satisfies $x'=Ax+\lambda F(x)$ which completes the proof.
\end{proof}

\section{Application to control theory}\label{control.theory}
Theorem \ref{theorem.Butcher.controle} provides a precise description of the solutions of \refeq{PF} when the data $u=(x^0,f)$ is 
a sum $u=\sum_{b\in\tree} \lambda^{\vert b\vert}u(b)$. We can use it in order to \emph{control} the solution of 
\refeq{PF} using $x^0$ or $f$ into the form of a sum indexed by planar trees. \\

Let us consider the following problem: given $m\in\N$, $B\in\mc{M}_{m,n}(\R)$ a $m\times n$ matrix and $x^0\in\R^n$, is there 
a function $v\in L^2((0,T),\R^m)$ such that the solution $x$ of 
\begin{equation}\label{PF2}\tag{$\mc{P}_\lambda$}
\left\{
\begin{array}{l}
\displaystyle{x\in\sol}\\
\displaystyle{x'=Ax+Bv+\lambda F(x)}\\
\displaystyle{x(0)=x^0\in\R^n}
\end{array}
\right.
\end{equation}
satisfies $x(T)=0$ ? \\

We show that if the corresponding linear system is controllable then we can define explicitly a control $v\in L^2((0,T),\R^m)$ 
as a sum indexed by planar trees $v=\sum_{b\in\tree}\lambda^{\vert b\vert} v(b)$ where for all $b\in\tree$, $v(b)$ is obtained by 
minimizing a functional $J(b)$. More precisely we have the following theorem 

\begin{theorem}\label{def.controle}\sl{
If $B$ satisfies the \emph{Kalman condition}\footnote{
This condition ensures \cite{Zuazua.Micu} that the problem can be solved if $\lambda=0$ \ie for all $x^0$ there is a function 
$v\in L^2((0,T),\R^m)$ such that the solution $x$ of \refeq{PF} with $\lambda=0$ corresponding to $x^0$ and $f=Bv$ satisfies 
$x(T)=0$. 
} 
\begin{equation}\label{Kalman}
\text{rank}\left( B,AB,\ldots,A^{n-1}B\right)=n
\end{equation}
Then there exists a family $(v(b))_{b\in \tree}$ of elements of $L^2((0,T),\R^m)$ such that for all $b\in\tree$, the function 
$v(b)$ is defined by $v(b):=B^*\tilde{y}(b)$, where $\tilde{y}(b)$ is the solution of the adjoint problem 
\begin{equation}\label{duale}\tag{$\mc{P}'$}
\left\{
\begin{array}{l}
\displaystyle{y\in\sol}\\
\displaystyle{-y'=A^*y}\\
\displaystyle{y(0)=y^0\in\R^n}
\end{array}
\right.
\end{equation}
corresponding to the initial condition $y^0=\tilde{y}^{(b)}\in\R^n$ which minimizes the functional 
$J(b):\R^n\longrightarrow\R$ defined by 
\begin{gather}
J(\circ)(y^0):=J_0(y^0)=\frac{1}{2}\int_0^T\vert B^* \tilde{y}(t)\vert^2\dd t +\lang x^0,y^0\rang\label{J0}\\
J(b)(y^0):=\frac{1}{2}\int_0^T\vert B^* \tilde{y}(t)\vert^2\dd t+
\int_0^T\lang \tilde{y}(t),F\left[(\Phi\ast u)(B_-(b))\right]\rang \dd t. \label{Jb}
\end{gather}
Here $\tilde{y}$ denotes the solution of \refeq{duale} with initial condition $y^0$ and $u(b)$ denotes the element of $\init$ 
defined by $u(\circ):=(x^0,B v(\circ))$ if $b=\circ$ and $u(b):=(0,Bv(b))$ if $b\neq\circ$. 
}
\end{theorem}

\begin{remark}
Notice that the right hand side of \refeq{Jb} only involve a function $v(c)$ such that $\vert c\vert<\vert b\vert$. Hence if 
the $v(c)$'s such that $\vert c\vert<\vert b\vert$ are known, one can compute $J(b)$. 
\end{remark}

\begin{proof}(of theorem \ref{def.controle})\\
The only thing we have to show is that the functionals $J(b)$ are well defined \ie that they admit a minimizer. 

It is well known (see \eg \cite{Zuazua.Micu}, \cite{Lee.Markus}) that if Kalman condition \refeq{Kalman} is satisfied then 
there exists $c_T$, which depends only on $T$, such that if $y$ is a solution of problem \refeq{duale},
\begin{equation}\label{app.observable}
\int_0^T\vert B^* y\vert^2\dd t\ge c_T\vert y^0\vert^2
\end{equation}
Hence we get 
$
J(\circ)(y^0)=\frac{1}{2}\int_0^T\vert B^* \tilde{y}(t)\vert^2\dd t +\lang x^0,y^0\rang
\ge \frac{c_T}{2}\vert y^0\vert^2-\vert x^0\vert \vert y^0\vert
\to\infty
$ 
when $\NO{y^0}\to\infty$. So $J(\circ)$ has a minimizer. \\

Let $N\in\N^*$ and suppose that the $v(c)$'s are well defined for all $c\in\tree$ such that $\vert c\vert<N$. 
Let $b\in\tree$ be such that $\vert b\vert=N$. Since $N\ge 1$, we get $b\neq\circ$. Hence 
there exists $r\ge 2$ and $(b_1,\ldots,b_r)\in\tree^2$ such that $B_-(b)=b_1\bullet \cdots\bullet b_r$. Then 
for all $i\in\lent 1,r\rent$, we have $\vert b_i\vert\le N-1$, so $F_r\left[(\Phi\ast u)(b_1\bullet \cdots\bullet b_r)\right]$ is well defined and 
\begin{align*}
J(b)(y^0)=&\frac{1}{2}\int_0^T\vert B^* \tilde{y}\vert^2\dd t
+\int_0^T\lang\tilde{y}(t),F\left[(\Phi\ast u)(B_-(b))\right]\rang\dd t\\
\ge& \frac{c_T}{2}\vert y^0\vert^2-\alpha \NO{F\left[(\Phi\ast u)(B_-(b))\right]}_{L^2((0,T),\R^n)} 
\vert y^0\vert \to +\infty\text{ when }\vert y^0\vert\to \infty
\end{align*}
where $\alpha$ denotes the norm of the linear map: $y^0\in\R^n\longmapsto$ solution $y\in\sol$ of \refeq{duale}. Then since $J(b)$ 
is continuous, $J(b)$ admits a minimizer $\tilde{y}^{(b)}\in\R^n$, which completes the proof. 
\end{proof}

\begin{theorem}\label{controle.simple}\sl{
Let $x^0\in\R^n$ and $B\in\mc{M}_{m,n}(\R)$ satisfies Kalman condition \refeq{Kalman}. Consider the family $(v(b))_{b\in\tree}$ 
of theorem \ref{def.controle}. 
If $\lambda$ is small enough then the sum $v:=\sum_{b\in\tree}\lambda^{\vert b\vert}v(b)$ makes sense in $\sol$ and 
provides a control for problem \refeq{PF2} \ie the solution $x$ of \refeq{PF2} corresponding to $(x^0,v)$ satisfies $x(T)=0$.   
}
\end{theorem}
In fact we have a more precise result: 
\begin{theorem}\label{theorem.controle.dim.finie}\sl{Suppose that $B$ satisfies the condition \refeq{Kalman}, then the family 
$(v(b))_{b\in\tree}$ of theorem \ref{def.controle} satisfies the following:
\begin{enumerate}
\item There exists constants $C$ and $C'$ which depend on $A$, $B$ and $T$, such that if  
$$
C'\vert \lambda\vert \NO{\Phi(\circ)}\NO{u(\circ)}^{-1}\vert F\vert\left(C\NO{u(\circ)} \right)<1,
$$
then the sum $\sum_{b\in\tree}\vert \lambda\vert^{\vert b\vert}\NO{v(b)}$ converges and 
the power series $\vert u\vert:=\sum_{b\in\tree}\vert\lambda\vert^{\vert b\vert}\NO{u(b)}$ satisfies 
$
\displaystyle{
\vert u\vert\le
\frac{16\NO{u(\circ)}}{1-C'\vert \lambda\vert \NO{\Phi(\circ)}\NO{u(\circ)}^{-1}\vert F\vert\left(C\NO{u(\circ)} \right)}
}
$.
\item Moreover, if we have 
$\displaystyle{\vert \lambda\vert \vert u\vert^{-1}\ \vert F\vert\left(16 \NO{\Phi(\circ)}\vert u\vert\right)<1}$,  
then the sum $v=\sum_{b\in\tree}\lambda^{\vert b\vert} v(b)$ converges and the solution $x$ of \refeq{PF} corresponding to 
$x^0$ and $f=Bv$ satisfies $x(T)=0$. 
\end{enumerate}
}
\end{theorem}

\begin{remark}
Notice that if Kalman condition \refeq{Kalman} is satisfied, theorem \ref{def.controle} ensures that we can always define 
the functions $v(b)$, but the sum $\sum_{b\in\tree} \lambda^{\vert b\vert}v(b)$ may not converge. But the first variation 
of the functional $J(b)$ of theorem \ref{def.controle} shows that identities \refeq{un.intro} and 
\refeq{deux.intro} of the introduction are satisfied, hence we can see the sum $\sum_{b\in\tree}\lambda^{\vert b\vert}v(b)$ as a 
"formal" control which is a  "real" control if $\lambda$ is small enough. 
\end{remark}

\begin{proof}(of theorem \ref{theorem.controle.dim.finie})\\
\textbf{Convergence of $u$. } 
Let us focus on the convergence of the power series $\sum_b \vert \lambda\vert^{\vert b\vert}\NO{v(b)}$. 

By definition, we know that $v(b)=B^* \tilde{y}(b)$ where $\tilde{y}(b)$ is the solution of \refeq{duale} corresponding to 
initial condition $\tilde{y}^{(b)}\in\R^n$ which minimize $J(b)$; hence computing the first variation of $J(\circ)$ and $J(b)$, 
we get that for all $y\in\sol$ solution of \refeq{duale}, we have 
\begin{gather}
\lang y^0,x^0\rang+\int_0^T \lang y(t),Bv(\circ)(t)\rang=0\label{app.v(circ)}\\
\forall b\neq \circ\  ;\ \int_0^T\lang y(t),Bv(b)(t)\rang+\int_0^T\lang y(t),F\left[(\Phi\ast u)(B_-(b))\right]\rang\dd t=0.
\label{app.v(b)}
\end{gather}
Taking $y=\tilde{y}(\circ)$ in \refeq{app.v(circ)} we get 
$\NO{v(\circ)}_{L^2}^2=-\lang\tilde{y}^{(\circ)},x^0\rang\le\vert x^0\vert\vert \tilde{y}^{(\circ)}\vert$. But since 
we assume that Kalman condition \refeq{Kalman} is satisfied, there exists $c_T>0$ such that \refeq{app.observable} occurs and 
we finally get 
\begin{equation}\label{maj.v(circ)}
\NO{v(\circ)}_{L^2}\le \frac{\vert x^0\vert}{\sqrt{c_T}}
\end{equation}

Now let $y=\tilde{y}(b)$ in \refeq{app.v(b)}, is leads to 
\begin{equation}\label{premiere.etape.controle}
\NO{v(b)}_{L^2}\le \frac{\alpha}{\sqrt{c_T}}\NO{F\left[(\Phi\ast u)(B_-(b))\right]}_{L^2}
\end{equation}
where $\alpha$ denotes the norm of the linear map: $y^0\in\R^n\longrightarrow y\in L^2((0,T),\R^n)\subset\sol$ solution of 
\refeq{duale}. Let us focus on $\NO{F\left[(\Phi\ast u)(B_-(b))\right]}_{L^2}$. \\

Starting from the definition\ref{coproduit} of $\cop$, it is easy to show that 
\begin{equation}\label{commute.B-.cop}
\cop\circ B_-=(id\otimes B_-)\circ\cop
\end{equation}
Let $b$ belong to $\tree$, $b\neq \circ$ then we denote by $r$ the number of children of the root of $b$. 
Using definition\ref{def.Phi.controle} of $\Phi\ast u$ and identity \refeq{commute.B-.cop}, we get 
$$
\NO{F\left[(\Phi\ast u)(B_-(b))\right]}\le \NO{F_r}
\sum_{\substack{b_{(1)},b_{(2)}\\ b_{(1)}\propto b_{(2)}=b}}
\NO{\Phi(B_-(b_{(2)})}\NO{u(b_{(1)})}
$$
which, together with lemma \ref{lemme.existence.controle}, leads to 
\begin{equation}\label{j.en.ai.besoin}
\NO{F\left[(\Phi\ast u)(B_-(b))\right]}\le\\
\sum_{\substack{b_{(1)} ,b_{(2)}\neq\circ\\ b_{(1)}\propto b_{(2)}=b}}
\NO{\Phi(\circ)}^{N(b_{(2)})-1}
\NO{u(b_{(1)})}
\NO{F_r}\prod_{j\in I(B_-(b_{(2)}))}\NO{F_{r_{B_-(b_{(2)})}(j)}}.
\end{equation}
But if $b_{(1)}\propto b_{(2)}=b$ where $b_{(2)}\neq \circ$, $r$ is the number of children of the root of $b_{(2)}$ too. 
Hence we have 
$$
\NO{F_r}\prod_{j\in I(B_-(b_{(2)}))}\NO{F_{r_{B_-(b_{(2)})}(j)}}=\prod_{j\in I(b_{(2)})}\NO{F_{r_{b_{(2)}}(j)}}
$$ 
So using inequality \refeq{j.en.ai.besoin} and \refeq{premiere.etape.controle} we finally get 
\begin{equation}\label{base.rec.convergence.controle}
\NO{v(b)}_{L^2}\le \frac{\alpha}{\sqrt{c_T}}
\sum_{\substack{b_{(1)} ,b_{(2)}\neq\circ\\ b_{(1)}\propto b_{(2)}=b}}
\NO{\Phi(\circ)}^{N(b_{(2)})-1}
\NO{u(b_{(1)})}
\prod_{j\in I(b_{(2)})}\NO{F_{r_{b_{(2)}}(j)}}.
\end{equation}
Then starting from this last inequality, we prove recursively the following lemma 
\begin{lemma}\label{lemme.maj.u}\sl{
For all $b\in\tree$ we have 
$$
\NO{v(\circ)}\le \frac{1}{\NO{B}}\beta^{\NO{b}-1}\NO{u(\circ)}^{\NO{b}}\left((1+\beta)\NO{\Phi(\circ)}\right)^{N(b)-1}
\prod_{j\in I(b)}\NO{F_{r_b(j)}}$$
where $\beta$ denotes the quantity 
$
\displaystyle{\beta:=\frac{\NO{B}\alpha}{\sqrt{c_T}}}
$. 
}
\end{lemma}
\begin{proof}(of lemma \ref{lemme.maj.u})\\
Let $\varphi$ denote the morphism of algebra (or the character) $\varphi:\foret\longrightarrow\R$ such that 
$\varphi(\circ):=\NO{u(\circ)}$ and for all $b\in\tree$, $b\neq \circ$, 
\begin{equation}\label{def.C(b)}
\varphi(b):=\beta\sum_{\substack{b_{(1)} ,b_{(2)}\neq\circ\\ b_{(1)}\propto b_{(2)}=b}}
\NO{\Phi(\circ)}^{\NO{b_{(2)}}-1}\varphi(b_{(1)}).
\end{equation}
Notice that if $b_{(2)}=\circ$ then $b_{(1)}=b$, so the definition of $\varphi$ shows directly that for all $b\in\tree$ 
\begin{equation}\label{def.C(b).2}
\varphi(b)=\frac{\beta}{1+\beta}\sum_{\substack{b_{(1)} ,b_{(2)}\\ b_{(1)}\propto b_{(2)}=b}}
\NO{\Phi(\circ)}^{\NO{b_{(2)}}-1}\varphi(b_{(1)}).
\end{equation}

Let us show recursively that for all $b\in\tree$, we have 
\begin{equation}\label{maj.u}
\NO{u(b)}\le \varphi(b)\NO{\Phi(\circ)}^{\vert b\vert }\prod_{j\in I(b)}\NO{F_{r_b(j)}}.
\end{equation}
If $b=\circ$ then definition $\varphi(\circ):=\NO{u(\circ)}$ shows that \refeq{maj.u} is satisfied. Let $b$ belong to 
$\tree$, $b\neq\circ$ and assume that \refeq{maj.u} is satisfied by all planar trees $c$ such that 
$\vert c\vert<\vert b\vert$. Since $b\neq\circ$, we have $u(b)=(0,Bv(b))$ so $\NO{u(b)}\le \NO{B}\NO{v(b)}$. Hence using 
\refeq{base.rec.convergence.controle}, we get 
\begin{equation}\label{trucmuche.fin}
\NO{u(b)}\le \beta
\sum_{\substack{b_{(1)} ,b_{(2)}\neq\circ\\ b_{(1)}\propto b_{(2)}=b}}
\NO{\Phi(\circ)}^{N(b_{(2)})-1}
\NO{u(b_{(1)})}
\prod_{j\in I(b_{(2)})}\NO{F_{r_{b_{(2)}}(j)}}.
\end{equation}
But all planar trees $c\in\tree$ which appear in $b_{(1)}$ on the right hand side of this last estimation 
satisfy $\vert c\vert<\vert b\vert$. So we finally get 
$$
\NO{u(b)}\le \beta
\sum_{\substack{b_{(1)} ,b_{(2)}\neq\circ\\ b_{(1)}\propto b_{(2)}=b}}
\NO{\Phi(\circ)}^{N(b_{(2)})-1}
\varphi(b_{(1)})\NO{\Phi(\circ)}^{\vert b_{(1)}\vert }\prod_{j\in I(b_{(1)})}\NO{F_{r_{b_{(1)}}(j)}}
\prod_{j\in I(b_{(2)})}\NO{F_{r_{b_{(2)}}(j)}}.
$$
If $b=b_{(1)}\propto b_{(2)}$ then the set $I(b)$ is composed of the internal vertices of $b_{(1)}$ and $b_{(2)}$. Moreover, 
we have $\vert b_{(1)}\vert+\vert b_{(2)}\vert=\vert b\vert$. Hence we get 
$$
\NO{u(b)}\le \beta\NO{\Phi(\circ)}^{\vert b\vert }\prod_{j\in I(b)}\NO{F_{r_b(j)}} 
\sum_{\substack{b_{(1)} ,b_{(2)}\neq\circ\\ b_{(1)}\propto b_{(2)}=b}}\NO{\Phi(\circ)}^{\NO{b_{(2)}}-1}\varphi(b_{(1)})
$$
where we recognize definition \refeq{def.C(b)} of $\varphi(b)$, so \refeq{maj.u} is true for $b$.\\ 

Let us study $\varphi(b)$. Again starting from the definition \refeq{coproduit} of $\cop$, we have 
\begin{equation}\label{lemme2}
\cop\circ B_+=(id\otimes B_+)\circ\cop+B_+\otimes \circ.
\end{equation}
Let $r$ belongs to $\R^*$ and $(b_1,\ldots,b_r)\in\tree^r$. Then using identity \refeq{lemme2} and definition 
\refeq{def.C(b)} of $\varphi$, we get 
$$
\varphi(B_+(b_1,\ldots,b_r))=\beta 
\sum_{\substack{b^1_{(1)} ,b^1_{(2)}\\ b^1_{(1)}\propto b^1_{(2)}=b_1}}\cdots
\sum_{\substack{b^r_{(1)} ,b^r_{(2)}\\ b^r_{(1)}\propto b^r_{(2)}=b_r}}
\NO{\Phi(\circ)}^{\NO{B_+(b^1_{(2)},\ldots,b^r_{(2)})}-1}\varphi(b^1_{(1)})\cdots \varphi(b^r_{(1)}).
$$
But since $\NO{B_+(b^1_{(2)},\ldots,b^r_{(2)})}=\NO{b^1_{(2)}}+\cdots+\NO{b^1_{(2)}}$, the last identity together with 
\refeq{def.C(b).2} leads to 
$$
\varphi(B_+(b_1,\ldots,b_r))=\beta\left(\frac{1+\beta}{\beta}\right)^r\NO{\Phi(\circ)}^{r-1} \prod_{j=1}^r \varphi(b_j).
$$
Then we can show very easily by recursion that $\varphi(b)$ is given by the following expression: 
\begin{equation}\label{C(b).fin}
\varphi(b)=\NO{u(\circ)}\left(\beta\NO{\Phi(\circ)}\NO{u(\circ)}\right)^{\NO{b}-1}(1+\beta)^{N(b)-1}
\end{equation}
and inserting this expression in \refeq{maj.u}, we finally get that for all $b\in\tree$ 
$$
\NO{u(\circ)}\le \beta^{\NO{b}-1}\NO{u(\circ)}^{\NO{b}}\left((1+\beta)\NO{\Phi(\circ)}\right)^{N(b)-1}
\prod_{j\in I(b)}\NO{F_{r_b(j)}}.
$$

Let us finish the proof: if $b=\circ$ then lemma \ref{lemme.maj.u} is obvious and if $b\neq\circ$ then 
\refeq{base.rec.convergence.controle} leads to 
$$
\NO{u(b)}\le \frac{\beta}{\NO{B}}
\sum_{\substack{b_{(1)} ,b_{(2)}\neq\circ\\ b_{(1)}\propto b_{(2)}=b}}
\NO{\Phi(\circ)}^{N(b_{(2)})-1}
\NO{u(b_{(1)})}
\prod_{j\in I(b_{(2)})}\NO{F_{r_{b_{(2)}}(j)}}.
$$
Then using \refeq{maj.u} and following the same steps, we get 
$$
\NO{v(b)}\le \frac{1}{\NO{B}} \varphi(b)\NO{\Phi(\circ)}^{\vert b\vert }\prod_{j\in I(b)}\NO{F_{r_b(j)}}
$$
which completes the proof of lemma \ref{lemme.maj.u} together with \refeq{C(b).fin}. 
\end{proof}\\

Now we can focus on the convergence of the sum $\sum_b\vert\lambda\vert^{\vert b\vert}\NO{v(b)}$. Let $N\in\N^*$,  
using lemma \ref{lemme.maj.u} and the fact that the number of planar trees $c$ such that $N(c)=N$ is bounded by $16^N$, we get 
$$
\sum_{\substack{b\in\tree(2,\infty)\\N(b)=N}}\vert\lambda\vert^{\vert b\vert}\NO{v(b)}\le 
\frac{16\NO{u(\circ)}}{\NO{B}} q_N
$$ 
where $q_N$ denotes the quantity 
$$
q_N:=\left(16(1+\beta)\beta\NO{u(\circ)}\right)^{N-1}\sum_{p=0}^N
\sum_{\substack{(r_1,\ldots,r_p)\in\N^p\\  r_1+\cdots+r_p=N-1}}
\left[\beta^{-1}\vert \lambda\vert\NO{\Phi(\circ)}\NO{u(\circ)}^{-1}\right]^p
\prod_{i=1}^p\NO{F_{r_i}}
$$
As in the proof of theorem \ref{theorem.Butcher.controle}, we consider the formal power series $Q:=\sum_{N\ge 1}q_N X^N$ 
and we get that 
\begin{equation}\label{fin}
Q=\frac{X}{1-\beta^{-1}\vert \lambda\vert \NO{\Phi(\circ)}\NO{u(\circ)}^{-1}\vert F\vert\left(
16(1+\beta)\beta\NO{u(\circ)} X
\right)}.
\end{equation}
Hence if 
$$
\beta^{-1}\vert \lambda\vert \NO{\Phi(\circ)}\NO{u(\circ)}^{-1}\vert F\vert\left(
16(1+\beta)\beta\NO{u(\circ)} 
\right)<1
$$
then the radius of convergence of $Q$ is greater than $1$, which shows that the sum 
$\sum_b\vert\lambda\vert^{\vert b\vert}\NO{v(b)}$ converges and moreover \refeq{fin} leads to  
$$
\vert v\vert:=\sum_{b\in\tree}\vert \lambda\vert^{\vert b\vert}\NO{v(b)}\le \sum_N v_N=
\frac{1}{\NO{B}}\frac{16\NO{u(\circ)}}{1-\beta^{-1}\vert \lambda\vert \NO{\Phi(\circ)}\NO{u(\circ)}^{-1}\vert F\vert\left(
16(1+\beta)\beta\NO{u(\circ)} X
\right)}.
$$
We get the estimation of $\vert u\vert$ by noticing that 
$\sum_{\substack{b\in\tree(2,\infty)\\N(b)=N}}\vert\lambda\vert^{\vert b\vert}\NO{v(b)}\le 16\NO{u(\circ)} q_N $. \\

\textbf{Verification of the control. }
Let us focus on the second part of the theorem. First of all theorem \ref{theorem.Butcher.controle} ensures that 
the sum $\sum_{b\in\tree}\vert\lambda\vert^{\vert b\vert}\lNO{(\Phi\ast u)(b)}$ converges and that 
$x=\sum_{b\in\tree}\lambda^{\vert b\vert}(\Phi\ast u)(b)$ is the solution of problem \refeq{PF2}. Let us show that 
we have $x(T)=0$. 

Let $y\in\sol$ be a solution of \refeq{duale} then since $x$ is a solution of \refeq{PF2}, we get 
\begin{equation}\label{app.base.controle.optimal}
\lang y(T),x(T)\rang-\lang y^0,x^0\rang=\int_0^T\dd t\ \frac{\dd \lang y(t),x(t)\rang }{\dd t}
=\int_0^T \lang y(t),Bv(t)\rang\dd t+\lambda\int_0^T \lang y(t),F(x(t))\rang\dd t.
\end{equation}
But we have seen (see proof of theorem \ref{theorem.Butcher.controle}) that $\lambda F(x)$ writes in $L^2((0,T),\R^n)$: 
$$
\lambda F(x)=\sum_{b\in \tree}\lambda^{\vert b\vert}F\left[(\Phi\ast u)(B_-(b))\right].
$$ 
Hence \refeq{app.base.controle.optimal} leads to 
$$
\lang y(T),x(T)\rang =\lang y^0,x^0\rang +\sum_{b\in\tree}
\lambda^{\vert b\vert}\left(\int_0^T \lang y,Bv(b)\rang +\int_0^T \lang y,F\left[(\Phi\ast u)(B_-(b))\right]\rang \right).
$$
But identities \refeq{app.v(circ)} and \refeq{app.v(b)} ensure that the right hand side of this last identity vanishes, so 
for all $y$ solution of \refeq{duale} we have $\lang y(T),x(T)\rang=0$. Since we can choose an arbitrary $y(T)\in\R^n$ this last 
property ensures that $x(T)=0$. 
\end{proof}

\section*{Acknowledgements} 
The author is very grateful to Sandrine Anthoine for careful reading of the manuscript and Fr\'ed\'eric H\'elein 
for helpful remarks and suggestions. 


\bibliographystyle{unsrt}

\end{document}